\documentclass[12pt]{article}
\usepackage{amssymb}

\usepackage{amsmath,amsthm,bbm}
\usepackage{version}
\usepackage{color}

\definecolor{dark_purple}{rgb}{0.4, 0.0, 0.4}
\definecolor{dark_green}{rgb}{0.0, 0.7, 0.0}

\usepackage{latexsym}
\usepackage{amssymb}
\usepackage{euscript}
\def\mcc{M\raise.5ex\hbox{c}C}
\def\mccarthy{M\raise.5ex\hbox{c}Carthy}

\def\K{{\cal K}}

\def\a{\alpha}
\def\l{\lambda}
\def\z{\zeta}

\def\vare{\varepsilon}

\let\i=\infty

\def\={\ = \ }

\def\C{\mathbb C}

\def\T{\mathbb T}
\def\D{\mathbb D}

\def\be{\setcounter{equation}{\value{theorem}} \begin{equation}}
\def\ee{\end{equation} \addtocounter{theorem}{1}}
\def\beq{\begin{eqnarray*}}
\def\eeq{\end{eqnarray*}}
\def\se{\setcounter{equation}{\value{theorem}}} 
\def\att{\addtocounter{theorem}{1}}
\def\vs{\vskip 5pt}

\def\bp{{\sc Proof: }}
\def\ep{{}{\hfill $\Box$} \vskip 5pt \par}

\def\oec{{}{\hspace*{\fill} $\lhd$} \vskip 5pt \par}
\def\bl{\begin{lemma}}
\def\el{\end{lemma}}
\def\bt{\begin{theorem}}
\def\et{\end{theorem}}
\def\bprop{\begin{prop}}
\def\eprop{\end{prop}}
\def\bd{\begin{definition}}
\def\ed{\end{definition}}
\def\br{\begin{remark}}
\def\er{\end{remark}}
\def\bexer{\begin{exercise}}
\def\eexer{\end{exercise}}
\def\bfig{\begin{figure}}
\def\efig{\end{figure}}

\newtheorem{theorem}{Theorem}[section]
\newtheorem{prop}[theorem]{Proposition}
\newtheorem{lemma}[theorem]{Lemma}

\newtheorem{question}[theorem]{Question}
\newtheorem{remark}[theorem]{Remark}
\newtheorem{defin}[theorem]{Definition}

\def\be{\begin{equation}}
\def\ee{\end{equation}}

\renewcommand\O{\Omega}

\newcommand\calv{{\mathcal V}}
\newcommand\U{{\mathcal U}}
\renewcommand\k{{\kappa}}
\renewcommand\o{{\omega}}
\newcommand\RR{{\mathcal R}}

\newcommand\car{Carath\'eodory\ }
\newtheorem{problem}[theorem]{Problem}

\title{Carath\'eodory sets in the tridisk}
\author{\L{}ukasz Kosi\'nski\\
Jagiellonian University\\
 Krakow
 \thanks{Partially supported by Sheng grant no. 2023/48/Q/ST1/00048 of the
National Science Center, Poland}
\and
John E. M\raise.5ex\hbox{c}Carthy
\thanks{Partially supported by National Science Foundation Grant DMS 
DMS 2054199}\\
Washington University in St. Louis
}

\date{2 May 2025}
\begin{document}

\bibliographystyle{plain}
\maketitle

{\sc Abstract:}
We characterize all algebraic subsets of the tridisk that are Carath\'eodory sets, that is the intrinsic
Carath\'eodory metric on the set equals the Carath\'eodory metric for the tridisk.
We show that such sets are either retracts, or are isomorphic to one particular exceptional set.

\section{Introduction}
\label{secintro}

%\subsection{Statement of Results}
%\vskip 10pt

Let $\D$ denote the open unit disk in the complex plane, and let 
 $\rho$  denote the pseudo-hyperbolic metric on $\D$. 
Let $\calv$ be any set with a complex structure, ie where it makes sense to speak of
holomorphic functions, and let $ {\rm Hol}(\calv, \D)$ denote the holomorphic functions that map into $\D$. 
Carath\'eodory had the idea that these functions can be used to describe an intrinsic distance on $\calv$ \cite{car13}.
We define
\be
\label{eqa1}
c_\calv(\l_1, \l_2) \= \sup \{ \rho(\phi(\l_1), \phi(\l_2)) : \phi \in {\rm Hol}(\calv, \D) \} .
\ee
The distance is called the \car pseudo-metric on $\calv$; it is a metric if the
bounded holomorphic functions separate points (eg if $\calv$ is an analytic subvariety of a bounded domain).
If $\calv$ is a subvariety of a domain $\O$, it is immediate that $c_\calv \ge c_\Omega|_\calv$, since one is taking the
supremum over an a priori larger set of functions.
It is natural to ask under what conditions the intrinsic \car metric on $\calv$ agrees with the one inherited from
$\O$; this is tantamount to asking if the functions where the extremal is attained in \eqref{eqa1} extend to
holomorphic maps from all of $\Omega$ into $\D$.

\begin{defin}
%Let $(\O, \calv)$ be a Cartan pair. 
Let $\O$ be pseudoconvex, and $\calv$ an analytic subvariety, We say $\calv$ is a Carath\'eodory set for $\O$ if, for every 
pair $\l_1, \l_2$ in $\calv$, we have
\[
c_\calv(\l_1, \l_2) \= c_\O (\l_1, \l_2) .
\]
\end{defin}
We give a definition of analytic subvariety in Definition \ref{defb1}.
The definition of a Carath\'eodory set is due to the first author and W. Zwonek in
\cite{kz21}, where they showed that the set $\K$ below is a Carath\'eodory set for $\D^3$.
\be
\label{eqa2}
\K \= \{ (x,y,z) \in \D^3: x + y + z = xy + yz + zx \} .
\ee

Let $\Omega$ be a domain in $\C^d$. A subset $\calv \subseteq \O$ is called a {\em retract} of $\O$
if there is a holomorphic map $r: \O \to \O$ such that $r(\O) = \calv$ and $r |_\calv = {\rm id}|_\calv$.
If $\calv$ is a retract of $\Omega$, then it is automatically a \car set, since any candidate function $\phi$ in \eqref{eqa1} extends to the function $\phi \circ r \in {\rm Hol}(\O,\D) \}$.

The main result of this note is that if $\calv$ is a Carath\'eodory set for $\D^3$ and $\calv$ is algebraic, then these are the only two possibilities: either $\calv$ is a retract or it is
biholomorphic to $\K$. By an algebraic subset of $\D^3$ we mean the intersection of $\D^3$ with an algebraic subset of $\C^3$.
\begin{theorem}
\label{thma1}
Let $\calv$ be an algebraic subset of $\D^3$. If $\calv$ is a Carath\'eodory set for $\D^3$, then either
$\calv$ is a retract of $\D^3$, or, up to a biholomorphism of $\D^3$, $\calv = \K$.
\end{theorem}

We prove this theorem in Section \ref{secc}. 

In \cite{kz21}, the authors proved that $\K$ is a totally geodesic subset of $\D^3$, 
which means that if $\phi : \D \to \D^3$ is a holomorphic map with a left inverse,  and $\phi(\D)$ is tangent to $\K$ at one point,
then $\phi(\D) \subset \K$. Moreover  $\K$ is not a graph over a product of lower dimensional totally geodesic sets, so it is in some sense primitive. See \cite{mmw17}
for another primitive two dimensional set, in the context of Teichm\"uller theory.
Theorem \ref{thma1} is further evidence that $\K$ is a special set.

%Let $\Omega$ be a domain in $\C^d$. A subset $\calv \subseteq \O$ is called a {\em retract} of $\O$
%if there is a holomorphic map $r: \O \to \O$ such that $r(\O) = \calv$ and $r |_\calv = {\rm id}|_\calv$.
%If $\calv$ is a retract of $\Omega$, then every holomorphic map $\phi : \calv \to \Sigma$, where $\Sigma$ is any
%set in $\C^n$, has an extension to a holomorphic map $\Phi: \O \to \Sigma$, namely
%\[
%\Phi \= \phi \circ r .
%\]
%Moreover the extension operator $E : \phi \mapsto \phi \circ r$ is automatically linear, and multiplicative
%if $\Sigma$ is equal to the unit disk $\D$.
%To what extent do properties of extension operators characterize retracts?

\begin{defin}
\label{defa1}
Let $\calv$ be an analytic subvariety of a pseudoconvex domain $\O \subseteq \C^d$. 
We say that the pair $(\O,\calv)$ has the {\em extension property} if every holomorphic 
$\phi: \calv \to \D$ has an extension to a holomorphic map $\Phi : \O \to \D$.
%When $\Sigma$ is the unit disk $\D$, we will refer to the $\D$-extension property simply as 
%the {\em  extension property.}
\end{defin}
%
%The definition of analytic subvariety is given in Section \ref{secdef}.
%We are interested in understanding what the $\Sigma$-extension property
%says about $\calv$, and in particular when it means that
%$\calv$ must be a retract of $\O$.
%The choice of $\Sigma$ is very important. 
%At one extreme, if $\Sigma = \C$, then a famous theorem of H. Cartan \cite{car51}  says
%that whenever $\calv$ is an  analytic subvariety of a pseudoconvex set $\O$, then every holomorphic function on $\calv$ extends to a holomorphic function on $\O$; in other words every pair $(\O,\calv)$ has the $\C$-extension property.
%At the other extreme, it is a tautology that if $(\O,\calv)$ has the $\calv$-extension property, then $\calv$ is a retract. 

It is immediate that 

(i) $\calv$ is a retract  of $\O$ $\Rightarrow$ 

(ii) $(\O,\calv)$ has the extension property $\Rightarrow$

 (iii) $\calv$ is a \car set for $\O$.

Asking when these implications can be reversed leads to interesting geometric questions.
There are currently no known examples of \car sets that are known not to have the extension property.
It is an open question whether $(\D^3,\K)$ has the extension property.
Our second result is that if one requires the extension to be linear, then in many cases this forces $\calv$ to be a retract.
We prove the following result in Section \ref{secd}.
\bt
\label{thma2}
Let $\O$ be a bounded convex  balanced  domain, and $\calv$ an analytic subvariety.
%Let $(\O, \calv)$ be a Cartan pair, and assume that $\O$ is bounded and convex.
Assume that  $\calv$ be relatively polynomially convex, and  that there is an isometric
linear extension operator $E: H^\i(\calv) \to H^\i(\O)$.
If any of the following three conditions hold, then $\calv$ is a retract of $\O$.

(i) $\O$ is the polydisk.

(ii) $\O$ is strictly convex.

(iii) $E$ is multiplicative.
\et

\section{History}

We say that $(\O, \calv)$ is a {\em Cartan pair} if $\O$ is a pseudoconvex domain, and $\calv$ is an analytic subvariety of $\O$.
By a holomorphic function on $\calv$, we mean a function $\phi$ with the property that for every $\l \in \calv$ there 
exists $\vare > 0$ so that $B(\l,\vare)$, the open ball centered at $\l$ with radius $\vare$, lies inside $\Omega$,
and there exists 
 a holomorphic function $\Phi$ on $B(\l,\vare)$ such that $\Phi |_{B(\l,\vare) \cap \calv}=  \phi_{B(\l,\vare) \cap \calv}$.
 It is far from obvious that a holomorphic function on $\calv$ extends to a holomorphic function on a neighborhood of $\calv$ in $\O$. Nevertheless, 
H. Cartan \cite{car51} proved the stronger result that if $(\O,\calv)$ is a Cartan pair, then any holomorphic function on $\calv$ extends to a holomorphic function on $\O$, though the extension may have a strictly larger range.

 W.  Rudin \cite[7.5.5]{rud69}
proved that if  $\O$ is the polydisk $\D^d$, and $\calv$ is an embedded polydisk of lower dimension inside $\O$, then if there is an extension operator $E: H^\i (\calv) \to H^\i(\O)$ that is either linear and isometric, or linear and multiplicative, 
 then $\calv$ is a retract of $\O$.

In \cite{agmcvn}, this result was strengthened for the bidisk.
\bt \cite{agmcvn} Let $\calv$ be a nonempty relatively polynomially convex subset of $\D^2$.
The  pair $(\D^2,\calv)$ has the  extension property if and only if $\calv$ is a  retract of $\D^2$.
\label{thmcit1}
\et

 In \cite{ghw08}, K. Guo, H. Huang and K. Wang improved Rudin's result for all $d$.
\bt
\cite{ghw08}
\label{thmcit2}
Suppose $\calv \subseteq \D^d$.

(i) If $d=3$, $\calv$ is an algebraic set,  $(\D^d,\calv)$ has the  extension property, and the extension can be chosen to be linear, then $\calv$ is a retract of $\D^3$.

(ii) For any $n$, if $\calv$ is $H^\i(\D^d)$ convex and there is an extension operator $E : H^\i (\calv) \to H^\i(\O)$ that is both linear and multiplicative (but not a priori required to be isometric),
then $\calv$ is a retract of $\D^d$.
\et

In \cite{maci19}, K. Maciaszek settled the issue for all one dimensional algebraic subsets of the polydisk.
\bt \cite{maci19}
\label{thmkm19}
If $\calv$ is a one dimensional algebraic subset of $\D^d$, 
and $(\D^d,\calv)$ has the  extension property, then $\calv$ is a retract. 
\et

 In \cite{kmc19}, the following was proved. 
 \bt
 \cite{kmc19}
 \label{thmcit3}
 Suppose $(\O,\calv)$ is a Cartan pair, and $\calv$ is relatively polynomially convex.
 If $(\O,\calv)$ has the extension property, then $\calv$ is a retract if any of the following 
 extra conditions hold:
 
 (i) $\O$ is the ball in any dimension.
 
 (ii) $\O$ is a strictly convex bounded subset of $\C^2$.
 
 (iii) $\O$ is a strongly linearly convex bounded subset of $\C^2$ with $C^3$ boundary. 
 \et

In \cite{maci20}, K. Maciaszek considered the Cartan domain ${\mathcal R}_{II}$ of
symmetric contractive $2$-by-$2$ matrices.
\bt
\label{thmcit4}
\cite{maci20}
Let $\O = \{ A \in {\mathcal M}_{2 \times 2} (\C) : A = A^t, I - A^* A \geq 0 \}$.
Let $\calv$ be an algebraic subset of $\O$. 
If $(\O,\calv)$ has the  extension property, then $\calv$ is a retract of $\O$.
\et

If one examines the proofs of Theorems \ref{thmcit1}, \ref{thmkm19}, \ref{thmcit3}, \ref{thmcit4}, they all remain true if ``$\calv$ has the extension property'' is replaced by the weaker ``$\calv$ is a \car set''.

The extension property does not always force $\calv$ to be a retract.
J. Agler, Z. Lykova and N. Young characterized subsets of the symmetrized bidisk with the extension property \cite{aly19}.
 \bt
 \cite{aly19}
 Let $\O = \{ (z+w, zw) : z,w \in \D \}$, and let 
 $\calv$ be an algebraic subset. Then $(\O,\calv)$ has the extension property 
   if and only if either
 $\calv$ is a retract of $\O$, or $\calv = {\mathcal R} \cup {\mathcal D}_\beta$, where
 ${\mathcal R} = \{(2z,z^2) : z \in \D\}$ and ${\mathcal D}_\beta = \{ (\beta + \bar \beta z, z):
 z \in \D \}$, and $\beta \in \D$.
 \et

In \cite{akmc22}, it was shown that there is no restriction whatsoever on $\calv$, if $\O$ is allowed to be arbitrary.
\bt
\cite{akmc22}
Let $(\O,\calv)$ be a Cartan pair.
Then there is a pseudoconvex domain $\Omega_1$ satisfying $\calv \subseteq \Omega_1 \subseteq \O$ 
so that $(\O_1,\calv)$ has the extension property.
\et
This result was extended in \cite{akmc23} to matrix-valued functions.

\section{Definitions and auxiliary results}
\label{secdef}

\begin{defin}
\label{defb1}
If $\Omega$ is a domain of holomorphy in $\C^d$, then an \emph{analytic subvariety of $\Omega$} is a relatively closed set $\calv$ in $\Omega$ such that for each point $\lambda\in \calv$ there exist a neighborhood $U$ of $\lambda$ ($U$ is open in $\C^d$) and holomorphic functions $f_1, f_2,\ldots, f_m$ on $U$ such that
\[
\calv\cap U =\{\mu\in U : f_i(\mu)=0 \text{ for } i=1,2,\ldots,m\}.
\]
An {\em algebraic subset of $\O$} is a set $W$ with the property that there exist 
polynomials $p_1, \dots, p_m$ so that 
\[
W \= \{ \l \in \O : p_j(\l) = 0, 1 \leq j \le m \} .
\]
\end{defin}

\begin{defin}
\label{defb2}
 Let $\calv \subseteq \C^d$. 
We say that a function $f:\calv\to \C$ is \emph{holomorphic on $\calv$} if for each $\lambda \in \calv$ there exist an open set $U \subseteq \C^d$ containing $\lambda$ and a holomorphic function $F$ defined on $U$ such that $F(\mu) =f(\mu)$ for all $\mu \in \calv \cap U$. We shall let ${\rm Hol}(\calv,\Sigma)$ denote the holomorphic maps on $\calv$ that take values in $\Sigma$.
\end{defin}

We shall let $H^\i (\calv)$ denote the bounded holomorphic functions from $\calv$ to $\C$, with the supremum norm.

\begin{defin}
If $(\O,\calv)$ is a Cartan pair, we say that $\calv$ is {\em relatively polynomially convex} in $\O$ if the intersection of the polynomially convex hull
of ${\calv}$ with $\O$ is $\calv$; in other words for every point $\l \in \O \setminus \calv$ there is a polynomial $p$
so that $|p (\l) | > \sup_{\mu \in \calv} |p(\mu)|$. We say $\calv$ is {\em $H^\infty(\O)$ convex} if for every
point $\l \in \O \setminus \calv$ there is a function $f \in H^\i(\O)$ 
so that $|f (\l) | > \sup_{\mu \in \calv} |f(\mu)|$.
\end{defin}

%\begin{defin}
%If $(\O, \calv)$ is a Cartan pair and $E: {\rm Hol}(\calv, \Sigma) \to {\rm Hol}(\O,\Sigma)$ is an extension operator
%(ie. $E(\phi) |_\calv = \phi$), then we say $E$ is linear if whenever $\phi_1, \phi_2 $ and $c_1 \phi_1 + c_2 \phi_2$ are in ${\rm Hol}(\calv, \Sigma)$, then $E(c_1 \phi_1 + c_2 \phi_2) = c_1 E(\phi_1) + c_2 E(\phi_2)$.
%If $\Sigma = \D$, we say $E$ is multiplicative if $E(\phi_1 \phi_2) = E(\phi_1) E(\phi_2)$.
%\end{defin}

A major obstacle to understanding when extension properties force sets to be retracts is that for most 
domains there is no explicit description of the retracts.
They are known for the ball and the polydisk.
The retracts of the ball are the intersections of affine sets with the ball \cite{su74}.
The retracts of the polydisk $\D^d$ were described by L. Heath and T. Suffridge \cite{hs81}.
They proved that, up to a permutation of coordinates, every $n$-dimensional retract of
$\D^d$ is of the form
\[
\{ (z_1, \dots, z_n, f_{1}(z_1, \dots, z_n), \dots f_{d-n}(z_1, \dots, z_n)) \}
\]
where $f_1, \dots, f_{d-n}$ are in ${\rm Hol}(\D^n, \D)$.

Let $(\a_1,\a_2,\a_3)$ be a triple of complex numbers, not all zero.
Let
\[
\K_{(\a_1,\a_2,\a_3)} \= \{ (x,y,z) \in \D^3 : \a_1 x + \a_2 y + \a_3 z =
\bar \a_1 yz + \bar \a_2 xz + \bar \a_3 xy \} .
\]
These analytic sets were studied in \cite{kz21}.
There are two distinct classes. If the triple $(|\a_1|,|\a_2|,|\a_3|)$ do not form the sides
of a triangle (which means that for some permutation of $\{1,2,3\}$ we have
$|\a_{i_3}| \geq |\a_{i_1}| + |\a_{i_2}|$), then $\K_{(\a_1,\a_2,\a_3)}$ is a retract, and biholomorphic to $\D^2$.
If the triple $(|\a_1|,|\a_2|,|\a_3|)$ does form the sides
of a triangle, then K. Maciaszek showed \cite{maci23} that there is a biholomorphism of $\D^3$
that maps $\K_{(\a_1,\a_2,\a_3)}$ onto $\K$ from  \eqref{eqa2}.

Let $m_a$ denote the Mobius map
\be
\label{eqb1}
m_a(\z) \= \frac{a - \z}{1 - \bar a \z} .
\ee

\section{Carath\'eodory sets for $\D^3$}
\label{secc}

The object of this section is to prove Theorem \ref{thma1}.

{\bf Theorem \ref{thma1}. }{\em
Let $\calv$ be an algebraic subset of $\D^3$. If $\calv$ is a Carath\'eodory set for $\D^3$, then either
$\calv$ is a retract of $\D^3$, or, up to a biholomorphism of $\D^3$, $\calv = \K$.
}
\vs

We start with the following theorems from \cite{kmc20}.
\bt \cite[Thm. 5.1]{kmc20}
\label{thmc1}
Let $\calv$ be a one dimensional algebraic subset of $\D^3$ and assume that $(\D^3,\calv)$ has the extension property.
Then $\calv$ is a retract of $\D^3$.
\et
\bt
\cite[Thm. 6.1]{kmc20}
\label{thmc2}
Let $\calv$ be a two dimensional relatively polynomially convex subset of $\D^3$
 and assume that $(\D^3,\calv)$ has the extension property.
Then either $\calv$ is a retract, or, for each $r = 1,2,3$,  there is a domain $U_r \subseteq \D^2$
and a holomorphic function $h_r : U_r \to \D$ so that
\se\att
\begin{align}\label{eqc1}
\calv  \=& \{ (x, y, h_3(x, y) ) : (x, y) \in U_3 \}\\\nonumber
\=& \{ (x, h_2(x, z) , z) : (x, z) \in U_2 \}\\\nonumber
\=& \{ ( h_1(y, z), y, z ) : (y, z) \in U_1 \} .
\end{align}
\et
It was observed in \cite[Rem. 13]{kz21} that the hypotheses of both theorems can be weakened from
assuming $(\D^3,\calv)$ has the extension property to merely requiring that $\calv$ is a Carath\'eodory set for $\D^3$. 
Moreover, the following result was also proved.
\bt
\label{thmc25} \cite[Cor. 10]{kz21}
Let $\calv$ be a Carath\'eodory set for $\D^2$ that is relatively polynomially convex.
Then $\calv$ is a retract of $\D^2$.
\et

 We shall explain how to prove the theorems with the weakened hypothesis. The key idea, originally due
 to P. Thomas \cite{th03}, is that Carath\'eodory extremals must map $\calv$ to a dense subset of $\D$.
A Carath\'eodory extremal for $\calv$ is a function $\phi$ for which the supremum in \eqref{eqa1} is attained.
 \bt
 \label{thmc26}
 \cite[Prop. 7]{kz21}
Let $\calv$ be a Carath\'eodory set for $\O$,   and let $\phi$ be a
Carath\'eodory extremal for a pair  $\l_1, \l_2$  in $\calv$.
Then $\phi(\calv)$ is dense in $\D$.
\et
 
Now if one examines the proof of Theorem \ref{thmc1} in \cite{kmc20}, one sees that 
the extension property is only used in two ways.
First, it is frequently used to argue that various Pick extremal functions must map
$\calv$ densely into $\D$. This argument can be replaced by 
Theorem \ref{thmc26}.
Secondly,  in \cite[Lemma 5.6]{kmc20},
the situation is reduced to when $\calv$ can be embedded in $\D^2$.
Once we reach that point, we can use Theorem \ref{thmc25}.

Likewise, the proof of Theorem \ref{thmc2} in \cite{kmc20}
only uses the extension property to conclude that
Pick extremal maps have dense range in $\D$, and Theorem \ref{thmc26} 
asserts that we only need to assume that $\calv$ is a Carath\'eodory set to
have this. 

\begin{defin}
A pair of distinct points $(\l,\mu)$ in $ \D^d$ is called $n$-balanced, for $1 \leq n \leq d$, if, 
for some permutation  $(i_1, \dots, i_d)$  of $(1,\dots, d)$, we have
\[
\rho(\l_{i_1}, \mu_{i_1})   = \dots = \rho(\l_{i_n}, \mu_{i_n}) \geq \rho(\l_{i_{n+1}}, \mu_{i_{n+1}}) \geq \dots \geq 
\rho(\l_{i_d}, \mu_{i_d}) .
\]
We shall say the pair is $n$-balanced w.r.t. the first $n$ coordinates if $(i_1, \dots, i_n) = (1, \dots, n)$.
\end{defin}

{\sc Proof of Theorem \ref{thma1}}.
Suppose $\calv$ is a Carath\'eodory set for $\D^3$ that is not a retract.
We want to prove it is biholomorphic to $\K$. By the remarks above, we can assume that
$\calv$ is two dimensional and  of the form \eqref{eqc1}.
Since it is codimension $1$, there is a square-free polynomial $P \in \C[x,y,z]$ so that $\calv$ is the intersection
of $\D^3$ with the zero set of $P$.

By applying a biholomorphism of $\D^3$ we can assume that $0 \in \calv$. So we can assume that
there is a domain $\U \subseteq \D^2$ and a holomorphic function $\k : \U \to \D$ so that $\k(0,0) = 0$
and such that
\be
\label{eqc2}
\notag
\calv \= \{ (x,y,\k(x,y)) : (x,y) \in \U \} .
\ee
%Moreover, since $\calv$ is assumed to be a semi-algebraic set, its projection $\U$ onto the first two coordinates 
%is also semi-algebraic, by the Tarski-Seidenberg theorem.

Our strategy is to prove that $\k$ is actually rational of degree $1$ in both variables.
Let us write $\T$ for the unit circle.

\begin{lemma}
\label{lemc1}
Suppose that $(0,\l)$ is a $3$-balanced pair of points in $\calv$.
Then there are unimodular constants $\omega_2, \omega_3$ so that the disk
$\{ (\zeta, \omega_2 \zeta, \omega_3 \zeta) : \zeta \in \D \}$ lies in $\calv$.
\end{lemma}
\bp
By hypothesis, $|\l_1| = |\l_2| = |\l_3|$.
Choose  $\omega_2, \omega_3$ so that 
$\l_1 =  \overline{\omega_2} \l_2 =  \overline{\omega_3} \l_3$.
Define $\phi: \D^3 \to \D$ by 
\[
\phi(x,y,z) \= \frac 13 ( x + \overline{\omega_2} y + \overline{\omega_3} z) .
\]
As \[
\rho(\phi(0), \phi(\lambda)) \ = \  c_{\D^3} (0, \l) \= c_{\calv} (0, \l),
\]
we see that $\phi$ is a Carath\'eodory extremal.
By Theorem \ref{thmc26}, $\overline{\phi(\calv)} = \phi(\overline{\calv})$ contains $\overline{\D}$,
and therefore \[
\{ (\zeta, \omega_2 \zeta, \omega_3 \zeta) : \zeta \in \T \} \subseteq \overline{\calv} .
\]
As $\calv$ is relatively polynomially convex (indeed, algebraic), the conclusion of the lemma follows.
\ep
\begin{lemma}
\label{lemc2}
(i)
Let $\omega$ be a unimodular number so that 
\be
\label{eqc3}
| \k_x(0,0) + \omega \k_y(0,0) | \ \leq \ 1 .
\ee
Then $\{ (\zeta, \omega \zeta) : \zeta \in \D \} \subseteq \U$.

(ii) If there is equality in \eqref{eqc3}, then 
\[
\{ (\zeta, \omega \zeta, (\k_x(0,0) + \omega \k_y(0,0) \zeta) : \zeta \in \D \} \subseteq \calv .
\]

(iii) If $|\k_x(0,0)|+ |\k_y(0,0) | \leq 1$, then $\calv$ is a retract of $\D^3$.
\end{lemma}
\bp
Assertion (ii) is an infinitesimal version of Lemma \ref{lemc1}, and is proved the same way.
So assume we have a strict inequality in \eqref{eqc3}.
Consider the map $\psi (\zeta) = \k(\z,\omega \z)$, with domain $\{ \z : (\z,\omega \z) \in \U\}$.
If there is some point $\z_0 \neq 0$ with $|\psi(\z_0)| = |\z_0|$, then by Lemma \ref{lemc1}
there is a disk $\{ (\z, \omega \z, \omega_3 \z) : \z \in \D\}$ in $\calv$, and hence (i) follows.
Otherwise, $| \psi(\z) | < |\zeta|$ for all $\zeta \neq 0$, and now the conclusion of (i)  follows from the fact
that $\calv$ is closed in $\D^3$.
 Indeed,
let $E = \{ \zeta \in \D: (\z,\omega \z) \in \U\}$. Suppose that $\zeta_0 \in \partial E \setminus \partial \D$.
If $\zeta_n$ are points in $E$ that converge to $\zeta_0$, then $(\zeta_n, \omega \zeta_n , \kappa(\zeta_n, \omega \zeta_n))$
are points in $\calv$ that tend to a point in $\overline{\calv} \setminus \calv $ which is in $\partial \D^3$.
But since $|\zeta_0| < 1$, this cannot happen.

(iii) It follows from (i) and (ii) that $(\zeta, \omega \zeta)$ is in $\U$ for any $\zeta \in \D, \omega \in \T$.
%Note that $\U$ is pseudoconvex. \blue There may be a better way to see that $\U$ is pseudoconvex.
%\black Indeed,
%since $\calv$ is the intersection of an algebraic set with $\D^3$, 
% every point in $\partial \U$ is either in $\partial (\D^2)$
% %or is a point where $\kappa$ cannot be analytically continued
% or is a point where $|\k|$ is $1$. So every boundary point is a point where some holomorphic function cannot be analytically continued, and hence $\U$ is pseudoconvex.
%\blue
%Is this valid?
%\black
%
%Let ${\mathcal D}$ be a Reinhardt domain that is contained in $\U$ and contains 
%$\cup \{ (\zeta, \omega \zeta) : \zeta \in \D, \omega \in \T \}$. 
%\blue
%Why does ${\mathcal D}$ exist? 
%\black
%Since $\U$ is pseudoconvex, it contains the
%envelope of holomorphy of ${\mathcal D}$, which, by \cite[Thm.1.12.4]{jp08}, contains $\D^2$.
%So we can conclude that $\U = \D^2$, and hence $\calv$ is a retract.
Therefore for any $0 < r < 1$ the power series for $\kappa$ converges on $r \D \times r \D$ to a function
of modulus less than $1$. So $\kappa$ extends to a holomorphic function from $\D^2$ to $\D$.
Since $\calv$ is the intersection of an algebraic set with $\D^3$, 
 every point in $\partial \U$ is either in $\partial (\D^2)$
 or is a point where $|\k|$ is $1$.
It follows that $\U = \D^2$, and $\calv$ is a retract.

\ep

Let the tangent plane to $\calv$ at $0$ be $\{ ax + by + c z = 0 \}$.
 The triple $(|a|,|b|,|c|)$ form the
sides of a triangle, since otherwise, by Lemma \ref{lemc2} (iii) applied to each pair of coordinates,
we would conclude that $\calv$ was a retract.
Let
\[
T \ := \ \{ \omega \in \T:  | \k_x(0,0) + \omega \k_y(0,0) |  < 1 \} .
\]
As $ \k_x(0,0) = - \frac ac$ and  $ \k_y(0,0) = - \frac bc$, the inequality $|b| + |c| > |a|$
implies that $T$ is an arc of positive length.
Note that if $\omega \in T$, then by Lemma \ref{lemc2} (i), the disk
$\{ (\zeta, \omega \zeta: \zeta \in \D \}$ is in $\U$. Hence the function $\zeta \mapsto \k(\z,\omega \z)$
is not only bounded by $1$ on $\D$, but actually extends to be continuous on $\overline{\D}$
since $\calv$ is algebraic.

Let
\[
W \ := \ 
\{ (x,y) \in \U : P_z(x,y,\k(x,y)) = 0 \} .
\]
Let $S$ be the subset of $T$ 
\[
S \ := \ \{ \omega \in T : \{ \xi \in \T, (\xi, \omega \xi) \in \overline{W} \} \ {\rm has\ positive\ measure}\}.
\]

{\sc Claim 1.} $S$ is finite.

Indeed, let $\omega \in S$, and define $f(\z) = P_z(\z, \omega \z, \k(\z, \omega \z))$.
Then $f$ is in the disk algebra, and vanishes on a set of positive measure on the boundary,
so must be identically zero. Therefore $\{ (\z, \o \z) : \z \in \D \} \subseteq W$.
But $W$ is  at most  one dimensional, since $P$ is square-free; 
and $W$ is semi-algebraic, by 
the Tarski-Seidenberg theorem, so it can only contain finitely many disks through $0$.
Therefore $S$ is finite, as claimed.
\oec

{\sc Claim 2.} If $\o \in T, \xi \in \T$ are such that $(\xi, \o \xi) \notin \overline{W}$,
then $|\k(\z,\o \z) | \to 1$ as $\z \to \xi$ from inside $\D$.

We prove this by contradiction. Suppose $\z_n$ converges to $\xi$ and 
$\k_n := \k(\z_n, \o \z_n)$ stays in a compact subset of $\D$.
Define
\[
k_n(x,y) \= m_{\k_n} (\k(m_{\z_n}(x), m_{\o \z_n}(y))) .
\]
Let 
\[
\calv_n = \{ (x,y,z): z = k_n(x,y),\quad (m_{\z_n}(x), m_{\o \z_n}(y)) \in \U \} .
\]
Then $\calv_n$ is the pushforward of $\calv$ by the automorphism of $\D^3$
given by $(m_{\z_n}, m_{\o \z_n} , m_{\k_n})$.
Since the points $\k_n$ remain in a compact set, we get
\be
\label{eqc4}
(k_n)_x(0,0) \ \approx \ (1 - | \z_n|^2) \k_x( \z_n, \o \z_n) .
\ee
Differentiating the equation $P(x,y, \k(x,y)) = 0$, we get
\be
\label{eqc5}
P_x( x,y, \k(x,y)) + P_z(x,y,\k(x,y)) \k_x(x,y) \= 0 .
\ee
The first term in \eqref{eqc5} remains bounded, and since 
$(\xi, \o \xi) \notin \overline{W}$, we have
$P_z( \z_n, \o \z_n, \k_n)$ does not tend to $0$.
Therefore there exists some constant $C$ so that 
\[
|\k_x(\z_n, \o \z_n) | \ \leq \ C \qquad \forall n .
\]
Plugging this into \eqref{eqc4}, we conclude that we can make $(k_n)_x(0,0)$ arbitrarily small
by choosing $n$ large enough. Interchanging $x$ and $y$, we can also make  $(k_n)_y(0,0)$
arbitrarily small.  Therefore by Lemma \ref{lemc2}, for $n$ large, we have that $\calv_n$ is a retract,
and hence so is $\calv$. This is our desired contradiction.
\oec

It follows that if $\o \in T \setminus S$, the function $\k(\z,\o\z)$ is in the disk algebra, and has modulus
one almost everywhere on the boundary, and hence is a finite Blaschke product. Let $n_\o$ be its degree.

Let $T_1$ be an arc in $T \setminus S$.
Since $\calv$ is algebraic, 
the function $\k(\xi, \o \xi)$ is continuous in $\omega$ for each fixed $\xi$ in $\T$,
so we can apply Rouch\'e's theorem to conclude that $n_\o$ is constant on $T_1$. 
By \cite[Thm. 5.2.2]{rud69}, this means that $\k(x,y)$ is a rational function of total degree $n = n_\o$.

{\sc Claim 3.} The function $\k$ is of degree $1$ in each variable. 

Indeed, we can also write $\calv$ as
\[
\calv \= \{ (h(y,z), y, z) : (y,z) \in \U_1 \}
\]
where $h$ is a rational function.
As the equality
\[
h(y, \k(x,y) ) \= x 
\]
holds for $(x,y) \in \U$, it extends to hold in $\C^2$.
It follows that $x \mapsto \k(x,y)$ is injective, and hence $\k$ is of degree $1$ in $x$.
Similarly, it is of degree $1$ in $y$.
\oec

{\sc Claim 4.} The total degree $n$ is $2$.

Otherwise, the degree is $1$. Since $\k(\z,\o\z)$ is a Blaschke product vanishing at $0$, we then
must have
\[
\k(x,y) \= a x + b y,
\]
for some $a,b$ satisfying $|a| + |b| > 1$ and $\left| |a| - |b| \right| < 1 $.
But in this case
 $T$ has positive length, and as 
$P(x,y,z) =C( z - ax - by)$ for some constant $C$, we get that $P_z = C$ is never $0$ and therefore
 $W$ is empty.
By Claim 2, we have that for every $\o$ in $T$ the function
\[
\k(\z,\o \z) = (1+\o) \z
\]
is inner. But we cannot have $|1+\o| =1$ for more than two values of $\o$,
so we have a contradiction.
\oec

Since $\k(0,0) = 0$, we can write
\be
\label{eqc55}
\k(x,y) \= \frac{a_1 x + a_2 y + a_3 xy}{1 + b_1 x + b_2 y + b_3 xy} .
\ee
As $k(\z,\o\z)$ is a Blaschke product of degree $2$ that vanishes at $0$, we get that
$|a_3| = 1, b_3 = 0, b_1 = \bar a_2 a_3, b_2 = \bar a_1 a_3$.
A straightforward calculation then shows that
\be
\label{eqc6}
\calv \= \{ (x,y,z) \in \D^3 : \a_1 x + \a_2 y + \a_3 z = \bar \a_1 yz + \bar \a_2 xz + \bar \a_3 xy \}
\ee
where, for some fixed choice $\beta$ of $\sqrt{\overline{a_3}}$, we have
$\a_1 = a_1 \beta, \a_2 = a_2 \beta, \a_3 = - \beta$.
It follows from \cite{maci23} that any such $\calv$ is biholomorphic to $\K$.
\ep

Remark: The argument in Claim 2 shows that if 
\[
\calv \= \{ (x,y,z) \in \D^3 : ax + by + cz = 0 \}
\]
and $\calv$ is a Carath\'eodory set 
then it is  a retract.
In \cite{kmc20} we used a 3 point argument to show that $z=x+y$ did not have the extension property.

\section{Linear extension operator from $H^\i(\calv)$ to $H^\i(\O)$}
\label{secd}

\begin{defin}
\label{defd1}
A domain $\O \subseteq \C^d$ is called strictly convex if every boundary point $\l$ is a point of strict convexity, i.e. there exists a hyperplane in $\C^d$ whose intersection with $\overline{\O}$ is $\{ \l \}$.
\end{defin}

In the proof of the following theorem, we shall only use the fact that we have an extension operator from polynomials restricted to $V$ to
bounded holomorphic functions on $\O$.
Recall that a set $\O$ is called balanced if $\l \O \subseteq \O$ for every $\l \in \overline{\D}$.
We need a minor variation of \cite[Thm. 3.7]{rud91}, which is stated for closed convex balanced sets;
a simple modification for open sets yields the following.

\begin{lemma}
\label{lemrud}
Let $\O$ be an open convex balanced set in a locally convex vector space $X$. Let $w \in X \setminus \O$.
Then there is a continuous linear functional $\Lambda$ on $X$ such that
$|\Lambda (z) | < 1 $ for all $z \in \O$, and $|\Lambda (w) | \ge 1$.
\end{lemma}
%\bp
%By \cite[Th. 3.4]{rud91} there exists $\Lambda \in X^*$ and $\gamma \in \R$ so that
%\[
%\Re \Lambda(z) \ <\  \gamma \ \le \  \Re \Lambda(w) \qquad \forall\ z \in \Omega .
%\]
%Since $\O$ is balanced, this means $|\Lambda(z)| < \gamma$ for all $z \in \Omega$.
%Rescaling, we can take $\gamma = 1$.
%\ep

{\bf Theorem \ref{thma2}.}
{\em
Let $(\O, \calv)$ be a Cartan pair, and assume that $\O$ is bounded, balanced,   and convex.
Let $\calv$ be relatively polynomially convex, and assume that there is an isometric
linear extension operator $E: H^\i(\calv) \to H^\i(\O)$.

If any of the following three conditions hold, then $\calv$ is a retract of $\O$.

(i) $\O$ is the polydisk.

(ii) $\O$ is strictly convex.

(iii) $E$ is multiplicative.
}

\vs
\bp
By the maximum principle, for any constant $c$, we must have $E(c) = c$.
%So without loss of generality, we can translate and assume $0 \in \calv$.
%Moreover,
Since $\calv$ is relatively polynomially convex and has the extension property, 
it is an analytic set in $\O$, by \cite[Lemma 2.1]{kmc20}.
For each coordinate function $z_j$, let $E(z_j) = \psi_j$. 
Let
$\psi = (\psi_1, \dots, \psi_d)$.
Note that $\psi(\O) \subseteq \O$.
Indeed, suppose $w $ were in the range of $\psi$ and not in $\O$.
By Lemma \ref{lemrud}  there would then be a linear functional $\Lambda$ on $\C^d$ that
would have modulus  that is less  than $1$ on $\overline{\O}$ but is greater than or equal to $1$ at $w$.
Writing $\Lambda(z) = \sum_{j=1}^d a_j z_j$ we would get that the function $\sum_{j=1}^d a_j \psi_j$ 
attains a value $w$ outside the  unit disk, whereas $ \sum_{j=1}^d a_j z_j$ maps $\calv$
into $\D$.
If  $|w| > 1$, this contradicts the isometric extension hypothesis. 
If $|w|=1$, it would violate the maximum modulus principle, since if $\sum_{j=1}^d a_j \psi_j$
were constant it would have to agree with its value on $\calv$, which has modulus less than one by hypothesis,
and if it is non-constant it cannot attain its maximum modulus at an interior point.

 Let  $\RR$ be the fixed point set of
$\psi$. By Vigu\'e's theorem \cite{vig85},
the fixed point set of a self-map of a bounded convex domain is a retract, so
 $\RR$ is a retract of $\O$. We have $\calv \subseteq \RR$;
we will show that in fact they are equal.

(i)  By  \cite{hs81}, every retract $\RR$ of $\D^d$  is,
after a permutation of coordinates, of the form
\[
\RR \= \{ (\zeta, f(\zeta) ) : \zeta \in \D^n \}
\]
for some holomorphic $f: \D^n \to \D^{d-n}$. Let $\pi_n$ denote projection onto the first $n$ coordinates. We claim that $\pi_n (\calv) = \D^n$.

Indeed, let $\tau = (\tau_1, \dots, \tau_n) \in \T^n$. Define
\[
L(z) = n + \sum_{k=1}^n \bar \tau_k z_k .\]
We have $E(L) = L \circ \psi$ which has norm $2n$ on $\D^d$. Therefore $L$ must have norm $2n$ on $\calv$, which means ${\rm cl}(\pi_n (\calv)) $ contains $\tau$. As $\tau$ is arbitrary,  ${\rm cl}(\pi_n (\calv)) $ contains $\T^n$. 

So $\calv$ is a holomorphic subvariety of $\D^d$ of the form
\[
\calv  \= \{ (\zeta, f(\zeta) ) : \zeta \in \U \}
\]
for some $\U \subseteq \D^n$ with ${\rm cl}(U) \supseteq \T^n$.

Suppose there exists some point
$\l = (\mu, f(\mu))  \in \RR \setminus \calv$. 
Since $\overline{\calv}$ is polynomially convex, there exists a polynomial $p$
with $| p(\l) | > \| p \|_\calv$.
 Let $\psi = E(p)$, and consider the function
\[
(\psi - p) \circ ({\rm id}_{\D^n} , f) : \D^n \to \C .
\]
This function vanishes on $\U$ and is non-zero at $\mu$, which contradicts  ${\rm cl}(\U) \supseteq \T^n$.

(ii)
We will show that $\overline{\RR} \cap \partial \O = \overline{\calv} \cap \partial \O$. Since $\overline{\calv}$ is polynomially convex,
this will mean $\overline{\RR} = \overline{\calv}$ and hence $\RR = \calv$.

Suppose there is some point $\l \in( \overline{\RR} \cap \partial \O) \setminus (\overline{\calv} \cap \partial \O)$.
Choose an affine polynomial $L$ that attains its maximum modulus $M$ on $\overline{\O}$ uniquely at $\l$.
Choose a sequence $\mu_n$ in $\RR$ with $\mu_n \to \l$.
Then $E(L) = L \circ \psi$, so 
\[
E(L)(\mu_n) = L(\psi(\mu_n)) = L(\mu_n) ,
\]
and $|L(\mu_n)| \to M$. This would mean that the norm of $E(L)$ is larger than $\| L \|_\calv$,
a contradiction.

(iii) Suppose there is some point $\l \in \RR \setminus \calv$. Choose a polynomial $p$ with
$|p(\l)| > \| p \|_\calv$. 
We have $E(p)(\l) = p \circ \psi (\l) = p(\l)$, which contradicts $E$ being isometric.
\ep

\section{Open Problems}

\begin{question}
\label{qf1}
Does $(\D^3, \K)$ have the extension property?
\end{question}

For the tridisk, exactly one of the converse implications (ii) $\Rightarrow$ (i) or (iii) $\Rightarrow$
(ii) after Definition \ref{defa1} holds. If the answer to Question \ref{qf1} is negative, 
then every algebraic set such that $(\D^3, \O)$ has the extension property is a retract.

\begin{question}
\label{qf2}
Does every function in ${\rm Hol}(\K, \D)$ extend to the Schur-Agler class of $\D^3$?
\end{question}

An affirmative answer to Question \ref{qf2} would give a variant of And\^o's inequality for triples of contractions satisfying the defining equation of $\K$.
It would clearly imply an affirmative answer to Question \ref{qf1}.

\begin{problem}
\label{qf3}
Construct a Cartan pair $(\O,\calv)$ that does not have the extension property but  so that $\calv$ is a Carath\'eodory set for $\O$.
\end{problem}

\vskip 10pt
{\bf Acknowledgement:} We would like to thank the referees for carefully reading the paper and suggesting improvements.

%\bibliography{../references}
%\bibliography{references}

\end{document}